\documentclass{hha}

\usepackage{amsmath}
\input diagxy
\newcommand{\defn}[1]{\emph{#1}}
\newcommand{\co}{\colon}
\newcommand{\categ}[1]{\mathsf{#1}}
\newcommand{\trunk}[1]{\mathsf{#1}}
\newcommand{\rmod}[1]{\mathcal{#1}}
\newcommand{\iso}{\cong}
\newcommand{\epic}{\twoheadrightarrow}

\newcommand{\Trunk}[1]{\operatorname{Trunk}(#1)}
\newcommand{\Ext}{\operatorname{Ext}}
\newcommand{\Hom}{\operatorname{Hom}}
\newcommand{\Conj}{\operatorname{Conj}}
\newcommand{\Sym}{\operatorname{Sym}}
\newcommand{\Id}{\operatorname{Id}}
\newcommand{\As}{\operatorname{As}}
\newcommand{\Op}{\operatorname{Op}}
\newcommand{\Orb}{\operatorname{Orb}}
\newcommand{\im}{\operatorname{im}}
\renewcommand{\epsilon}{\varepsilon}

\renewcommand{\theenumi}{\roman{enumi}}

\renewcommand{\to}{\rightarrow}

\newtheorem{example}{Example}[section]
\newtheorem{theorem}{Theorem}[section]
\newtheorem{corollary}[theorem]{Corollary}
\newtheorem{lemma}[theorem]{Lemma}
\newtheorem{proposition}[theorem]{Proposition}

\begin{document}

\title{Extensions of racks and quandles}

\author{Nicholas Jackson}

\email{nicholas@maths.warwick.ac.uk}
\address{Mathematics Institute\\
         University of Warwick\\
         Coventry\\
         CV4 7AL\\
         United Kingdom}
\classification{18G15; 18E10, 18G60}
\keywords{Racks, quandles, extensions, modules, homology, cohomology}
\begin{abstract}
A \defn{rack} is a set equipped with a bijective,
self-right-dist\-rib\-u\-tive
binary operation, and a \defn{quandle} is a rack which satisfies an
idempotency condition.

In this paper, we introduce a new definition of modules over a rack or
quandle, and show that this definition includes the one studied by Etingof
and Gra\~na \cite{etingof/grana:orc} and the more general one given by
Andruskiewitsch and Gra\~na \cite{andr/grana:pointed-hopf}.  We further
show that this definition coincides with the appropriate specialisation
of the definition developed by Beck \cite{beck:thesis}, and hence that
these objects form a suitable category of coefficient objects in which
to develop homology and cohomology theories for racks and quandles.

We then develop an Abelian extension theory for racks and quandles which
contains the variants developed by Carter, Elhamdadi, Kamada and
Saito \cite{carter/elhamdadi/saito:twisted,carter/kamada/saito:diag} as
special cases.
\end{abstract}

\received{January 19, 2005}   % receive date (for example: 11 October 1999)
\revised{May 18, 2005}    % receive date
\published{June 20, 2005}  % publish date
\submitted{Ronald Brown}  % Name of Journal's Editor, who submitted Article 

\volumeyear{2005} % Volume Year
\volumenumber{7}  % Volume Number 
\issuenumber{1}   % Issue Number

\startpage{151}     % PageNumber of first page

\maketitle

\section{Introduction}
A \defn{rack} (or \defn{wrack}) is a set $X$ equipped with a
self-right-distributive binary operation (often written as exponentiation)
satisfying the following two axioms:
\begin{enumerate}
\item[(R1)] For every $a,b \in X$ there is a unique $c \in X$
  such that $c^b = a$.
\item[(R2)] For every $a,b,c \in X$, the \defn{rack identity} holds:
$$a^{bc} = a^{cb^c}$$
\end{enumerate}
In the first of these axioms, the unique element $c$ is often denoted
$a^{\overline{b}}$, although $\overline{b}$ should not itself be regarded
as an element of the rack.  Association of exponents should be understood
to follow the usual conventions for exponential notation.  In particular,
the expressions $a^{bc}$ and $a^{cb^c}$ should be interpreted as $(a^b)^c$
and $(a^c)^{(b^c)}$ respectively.

A rack which, in addition, satisfies the following idempotency criterion
is said to be a \defn{quandle}.
\begin{enumerate}
\item[(Q)] For every $a \in X$, $a^a = a$.
\end{enumerate}

There is an obvious notion of a \defn{homomorphism} of racks: a
function $f\co X \to Y$ such that $f(a^b) = f(a)^{f(b)}$ for all $a,b \in
X$.  We may thus form the categories $\categ{Rack}$ and $\categ{Quandle}$.

For any element $x \in X$ the map $\pi_x\co a \mapsto a^x$ is a bijection.
The subgroup of $\Sym X$ generated by $\{ \pi_x : x \in X \}$ is the
\defn{operator group} of $X$, denoted $\Op X$.  This assignment is not
functorial since there is not generally a well-defined group homomorphism
$\Op f\co \Op X \to \Op Y$ corresponding to an arbitrary rack
homomorphism $f\co X \to Y$.  The group $\Op X$ acts on the rack $X$,
and divides it
into \defn{orbits}.  Two elements $x,y \in X$ are then said to be in
the same orbit (denoted $x \sim y$ or $x \in [y]$) if there is a (not
necessarily unique) word $w \in \Op X$ such that $y = x^w$.  A rack with
a single orbit is said to be \defn{transitive}.  The set of orbits of $X$
is denoted $\Orb X$.

Given any group $G$, we may form the \defn{conjugation rack} $\Conj G$
of $G$ by taking the underlying set of $G$ and defining the rack operation
to be conjugation within the group, so $g^h := h^{-1}gh$ for all $g,h
\in G$.  This process determines a functor $\Conj\co \categ{Group}
\to \categ{Rack}$ which has a left adjoint, the \defn{associated group}
functor $\As\co \categ{Rack} \to \categ{Group}$.  For a given rack $X$,
the associated group $\As X$ is the free group on the elements of $X$
modulo the relations
$$a^b = b^{-1}ab$$
for all $a,b \in X$.

Racks were first studied by Conway and Wraith \cite{conway/wraith:wracks}
and later (under the name `automorphic sets') by
Brieskorn \cite{brieskorn:automorphic}, while quandles were introduced
by Joyce \cite{joyce:knot-quandle}.  A detailed exposition may be found in the
paper by Fenn and Rourke \cite{fenn/rourke:racks-links}.

A \defn{trunk} $\trunk{T}$ is an object analogous to a category, and
consists of a class of \defn{objects} and, for each ordered pair $(A,B)$ of
objects, a set $\Hom_{\trunk{T}}(A,B)$ of \defn{morphisms}.  In
addition, $\trunk{T}$ has a number of \defn{preferred squares}
$$\bfig\square[A`B`C`D;f`g`h`k]\efig$$
of morphisms, a concept analogous to that of composition in a category.
Morphism composition need not be associative, although it is in all the
cases discussed in this paper, and particularly when the trunk in question
is also a category.

Given two arbitrary trunks $\trunk{S}$ and $\trunk{T}$, a \defn{trunk map} or
\defn{functor} $F\co\trunk{S}\to\trunk{T}$ is a map which assigns to
every object $A$ of $\trunk{S}$ an object $F(A)$ of $\trunk{T}$, and to
every morphism $f\co A \to B$ of $\trunk{S}$ a morphism $F(f)\co F(A) \to
F(B)$ of $\trunk{T}$ such that preferred squares are preserved:
$$\bfig\square/>`>`>`>/<500,500>[F(A)`F(B)`F(C)`F(D);f_*`g_*`h_*`k_*]\efig$$
For any category $\categ{C}$ there is a well-defined trunk
$\Trunk{\categ{C}}$ which has the same objects and morphisms as
$\categ{C}$, and whose preferred squares are the commutative diagrams
in $\categ{C}$.  In particular, we will consider the case
$\Trunk{\categ{Ab}}$, which we will denote $\trunk{Ab}$ where there is
no ambiguity.  Trunks were first introduced and studied by Fenn,
Rourke and Sanderson \cite{fenn/rourke/sanderson:trunks}.

In this paper, we study extensions of racks and quandles in more
generality than before, in the process describing a new, generalised
notion of a module over a rack or quandle, which is shown to coincide with
the general definition of a module devised by Beck \cite{beck:thesis}.
Abelian groups $\Ext(X,\rmod{A})$ and $\Ext_Q(X,\rmod{A})$ are defined
and shown to classify (respectively) Abelian rack and quandle extensions
and to be generalisations of all known existing $\Ext$ groups for racks
and quandles.

This paper contains part of my doctoral thesis \cite{jackson:thesis}.
I am grateful to my supervisor Colin Rourke, and to Alan Robinson,
Ronald Brown, and Simona Paoli for many interesting discussions and
much helpful advice over the past few years.  I also thank the referees for
their kind comments and helpful suggestions.

\section{Modules}
Given a rack $X$ we define a trunk $\trunk{T}(X)$ as follows:  let
$\trunk{T}(X)$ have one object for each element $x \in X$, and for each
ordered pair $(x,y)$ of elements of $X$, a morphism $\alpha_{x,y}\co x \to
x^y$ and a morphism $\beta_{y,x}\co y \to y^x$ such that the squares
$$\bfig\square/>`>`>`>/<1000,500>[x`x^y`x^z`x^{yz} = x^{zy^z};
  \alpha_{x,y}`\alpha_{x,z}`\alpha_{x^y,z}`\alpha_{x^z,y^z}]\efig
  \qquad
  \bfig\square/>`>`>`>/<1000,500>[y`x^y`y^z`x^{yz} = x^{zy^z};
  \beta_{y,x}`\alpha_{y,z}`\alpha_{x^y,z}`\beta_{y^z,x^z}]\efig$$
are preferred for all $x,y,z \in X$.

Thus a trunk map $A\co\trunk{T}(X)\to\trunk{Ab}$, as defined in the
previous section, determines Abelian groups $A_x$, and Abelian group
homomorphisms
$\phi_{x,y}\co A_x \to A_{x^y}$ and $\psi_{y,x}\co A_y \to A_{x^y}$,
such that
\begin{eqnarray*}
\phi_{x^y,z}\phi_{x,y} & = & \phi_{x^z,y^z}\phi_{x,z} \\
\mbox{and}\qquad
\phi_{x^y,z}\psi_{y,x} & = & \psi_{y^z,x^z}\phi_{y,z}
\end{eqnarray*}
for all $x,y,z \in X$.  It will occasionally be convenient to denote
such a trunk map by a triple $(A,\phi,\psi)$.

\subsection{Rack modules}
Let $X$ be an arbitrary rack. Then a \defn{rack module} over $X$ (or an
\defn{$X$--module}) is a trunk map
$\rmod{A} = (A,\phi,\psi)\co\trunk{T}(X) \to \trunk{Ab}$ such that each
$\phi_{x,y}\co A_x \iso A_{x^y}$ is an isomorphism, and
\begin{equation}
\label{eqn:rmod}
\psi_{z,x^y}(a) = \phi_{x^z,y^z}\psi_{z,x}(a) + \psi_{y^z,x^z}\psi_{z,y}(a)
\end{equation}
for all $a \in A_z$ and $x,y,z \in X$.

If $x,y$ lie in the same orbit of $X$ then this implies that $A_x \iso
A_y$ (although the isomorphism is not necessarily unique).  For racks with
more than one orbit it follows that if $x \not\sim y$ then $A_x$ need
not be isomorphic to $A_y$.  Rack modules where the constituent groups are
nevertheless all isomorphic are said to be \defn{homogeneous}, and those
where this is not the case are said to be \defn{heterogeneous}.  It is
clear that modules over transitive racks must be homogeneous.

An $X$--module $\rmod{A}$ of the form $(A,\Id,0)$ (so that $\phi_{x,y} =
\Id\co A_x \to A_{x^y}$ and $\psi_{y,x}$ is the zero map $A_y \to
A_{x^y}$) is said to be \defn{trivial}.

\begin{example}[Abelian groups]
\label{exm:abgroup}
{\rm Any Abelian group $A$ may be considered as a homogeneous trivial
$X$--module $\rmod{A}$, for any rack $X$, by setting
$A_x = A, \phi_{x,y} = \Id_A$, and $\psi_{y,x} = 0_A$ for all $x,y \in
X$.}
\end{example}

\begin{example}[$\As X$--modules]
\label{exm:asx}
{\rm Let $X$ be a rack, and let $A$ be an Abelian group equipped with an
action of $\As X$.  Then $A$ may be considered as a homogeneous
$X$--module $\rmod{A} = (A,\phi,\psi)$ by setting $A_x = A$, and defining
$\phi_{x,y}(a) = a \cdot x$ and $\psi_{y,x}(a) = 0$ for all $a \in A$
and $x,y \in X$.}
\end{example}

In particular, Etingof and Gra\~na \cite{etingof/grana:orc} study a
cohomology theory for racks, with $\As X$--modules as coefficient objects.

\begin{example}
\label{exm:andr/grana}
{\rm In \cite{andr/grana:pointed-hopf}, Andruskiewitsch and Gra\~na define an
\defn{$X$--module} to be an Abelian group $A$ equipped with a family $\eta
= \{\eta_{x,y} : x,y \in X\}$ of automorphisms of $A$ and another family
$\tau = \{\tau_{x,y} : x,y \in X\}$ of endomorphisms of $A$ such that
(after slight notational changes):
\begin{eqnarray*}
\eta_{x^y,z}\eta_{x,y} & = & \eta_{x^z,y^z}\eta_{x,z} \\
\eta_{x^y,z}\tau_{y,x} & = & \tau_{y^z,x^z}\eta_{y,z} \\
\tau_{z,x^y} & = & \eta_{x^z,y^z}\tau_{z,x} + \tau_{y^z,x^z}\tau_{z,y}
\end{eqnarray*}
This may readily be seen to be a homogeneous $X$--module in the
context of the current discussion.

As a concrete example, let $X$ be $C_3=\{0,1,2\}$, the cyclic rack with
three elements.  This has rack structure given by $x^y = x+1 \pmod{3}$
for all $x,y \in X$.  Let $A=\mathbb{Z}_5$ and define:
\begin{align*}
&\eta_{x,y}\co A \to A;\quad n \mapsto 2n\pmod{5}\\
&\tau_{y,x}\co A \to A;\quad n \mapsto 4n\pmod{5}
\end{align*}
Then this satisfies Andruskiewitsch and Gra\~na's definition of a
$C_3$--module, and (by setting $A_0=A_1=A_2=A=\mathbb{Z}_5$) is also a
homogeneous $C_3$--module in the context of the current discussion.}
\end{example}

\begin{example}[Alexander modules]
\label{exm:alexander}
{\rm Let $h = \{ h_i : i \in \Orb X \}$ be a family of Laurent polynomials
in one variable $t$, one for each orbit of the rack $X$, and let $n =
\{ n_i : i \in \Orb X \}$ be a set of positive integers, also one for
each orbit.  Then we may construct a (possibly heterogeneous) $X$--module
$\rmod{A} = (A,\phi,\psi)$ by setting $A_x =
\mathbb{Z}_{n_{[x]}}[t,t^{-1}]/h_{[x]}(t)$, $\phi_{x,y}\co a \mapsto ta$,
and $\psi_{y,x}\co b \mapsto (1-t)b$ for all $x,y \in X$, $a \in A_x$ and
$b \in A_y$.
The case where $A_x = \mathbb{Z}[t,t^{-1}]/h_{[x]}(t)$ for all
$x$ in some orbit(s) of $X$ is also an $X$--module.}
\end{example}

\begin{example}[Dihedral modules]
\label{exm:dihedral}
{\rm Let $n = \{ n_i : i \in \Orb X \}$ be a set of positive integers, one
for each orbit of $X$.  Then let $\rmod{D} = (D,\phi,\psi)$ denote the
(possibly heterogeneous) $X$--module where $D_x = \mathbb{Z}_{n_{[x]}}$,
$\phi_{x,y}(a) = -a$, and $\psi_{y,x}(b) = 2b$ for all $x,y \in X$, $a \in
A_x$ and $b \in A_y$.  This module is isomorphic to the Alexander module
where $h_i(t) = (1+t)$ for all $i \in \Orb X$.  The case where $A_x =
\mathbb{Z}$ for all $x$ in some orbit(s) of $X$, is also an $X$--module.
The \defn{$n$th homogeneous dihedral $X$--module} (where all the $n_i$ are
equal to $n$) is denoted $\rmod{D}_n$.  The case where $D_x = \mathbb{Z}$
for all $x \in X$ is the \defn{infinite homogeneous dihedral $X$--module}
$\rmod{D}_\infty$.}
\end{example}

Given two $X$--modules $\rmod{A} = (A,\phi,\psi)$ and $\rmod{B} =
(B,\chi,\omega)$, a \defn{homomorphism} of $X$--modules, or an
\defn{$X$--map}, is a natural transformation $f\co \rmod{A} \to \rmod{B}$
of trunk maps,
that is, a collection $f = \{ f_x\co A_x
\to B_x : x \in X \}$ of Abelian group homomorphisms such that
\begin{eqnarray*}
\phi_{x,y}f_x & = & f_{x^y}\phi_{x,y} \\
\mbox{and}\qquad\psi_{y,x}f_y & = & f_{x^y}\psi_{y,x}
\end{eqnarray*}
for all $x,y \in X$.

We may thus form the category $\categ{RMod}_X$ whose objects are
$X$--modules, and whose morphisms are $X$--maps.

In his doctoral thesis \cite{beck:thesis}, Beck gives a general
definition of a `module' in an arbitrary category.  Given a category
$\categ{C}$, and an object $X$ of $\categ{C}$, a \defn{Beck module} over $X$
is an Abelian group object in the slice category $\categ{C}/X$.
For any group $G$, the category $\categ{Ab}(\categ{Group}/G)$, for
example, is equivalent to the category of $G$--modules.  Similar results
hold for Lie algebras, associative algebras and commutative rings.  The
primary aim of this section is to demonstrate a categorical equivalence
between the rack modules just defined, and the Beck modules in the
category $\categ{Rack}$.

For an arbitrary rack $X$ and an $X$--module $\rmod{A} = (A,\phi,\psi)$,
we define the \defn{semidirect product} of $\rmod{A}$ and $X$ to be the set
$$\rmod{A} \rtimes X = \{ (a,x) : x \in X, a \in A_x \}$$
with rack operation given by
$$(a,x)^{(b,y)} := \left(\phi_{x,y}(a) + \psi_{y,x}(b), x^y\right).$$
\begin{proposition}
\label{thm:semidirect-rack}
\label{prp:semidirect-rack}
For any rack $X$ and $X$--module $\rmod{A} = (A,\phi,\psi)$, the
semidirect product $\rmod{A} \rtimes X$ is a rack.
\end{proposition}
\begin{proof}
For any three elements $(a,x), (b,y), (c,z) \in \rmod{A} \rtimes X$,
\begin{align*}
(a,x)^{(b,y)(c,z)}\! &= (\phi_{x,y}(a) + \psi_{y,x}(b),x^y)^{(c,z)} \\
&= (\phi_{x^y,z}\phi_{x,y}(a) + \phi_{x^y,z}\psi_{y,x}(b)
  + \psi_{z,x^y}(c),x^{yz}) \\
&= (\phi_{x^z,y^z}\phi_{x,z}(a) + \psi_{y^z,x^z}\phi_{y,z}(b)
  + \phi_{x^z,y^z}\psi_{z,x}(c) + \psi_{y^z,x^z}\psi_{z,y}(c), x^{zy^z}) \\
&= (\phi_{x,z}(a) + \psi_{z,x}(c),x^z)^{(\phi_{y,z}(b) + \psi_{z,y}(c),y^z)} \\
&= (a,x)^{(c,z)(b,y)^{(c,z)}}.
\end{align*}
Furthermore, for any two elements $(a,x), (b,y) \in \rmod{A} \rtimes X$,
there is a unique element
$$(c,z) = (a,x)^{\overline{(b,y)}} = (\phi_{z,y}^{-1}(a -
\psi_{y,z}(b)),x^{\overline{y}}) \in \rmod{A} \rtimes X$$
such that $(c,z)^{(b,y)} = (a,x)$.

Hence $\rmod{A} \rtimes X$ satisfies the rack axioms.
\end{proof}

\begin{theorem}
\label{thm:rmod-beck}
For any rack $X$, the category $\categ{RMod}_X$ of $X$--modules is
equivalent to the category
$\categ{Ab}(\categ{Rack}/X)$ of Abelian group objects over $X$.
\end{theorem}

\begin{proof}
Given an $X$-module $\rmod{A} = (A,\phi,\psi)$, let $T\rmod{A}$ be the
object $p\co \rmod{A} \rtimes X \epic X$ in the slice category
$\categ{Rack}/X$, where $p$ is defined as projection onto the second coordinate.
Given an $X$--map $f\co\rmod{A}\to\rmod{B}$, we obtain
a slice morphism
$Tf:\rmod{A} \rtimes X \to \rmod{B} \rtimes X$ defined by $T(f)(a,x) =
(f_x(a),x)$ for all $a \in A_x$ and $x \in X$.  This is functorial since,
for any $X$--module homomorphism $g\co\rmod{B}\to\rmod{C}$,
\begin{align*}
T(fg)(a,x) &= ((fg)_x(a),x) \\
&= (f_xg_x(a),x) \\
&= T(f)(g_x(a),x) \\
&= T(f)T(g)(a,x)
\end{align*}
for all $a \in A_x$ and $x \in X$.  We thus have a functor
$T\co\categ{RMod}_X\to\categ{Rack}/X$.  Our aim is to show firstly that
the image of $T$ is the subcategory $\categ{Ab}(\categ{Rack}/X)$, and
secondly that $T$ has a well-defined inverse.

To show the first, that $T\rmod{A}$ has a canonical structure as an
Abelian group object, we must construct an appropriate section, and
suitable multiplication and inverse morphisms.

Let:
\begin{align*}
& r\co \rmod{A} \rtimes X \to \rmod{A} \rtimes X  ; &&
     (a,x) \mapsto (-a,x) \\
& m\co (\rmod{A} \rtimes X)\times_X(\rmod{A} \rtimes X) \to
  \rmod{A} \rtimes X ; &&
     ((a_1,x),(a_2,x)) \mapsto (a_1+a_2,x) \\
& s\co X \to \rmod{A} \rtimes X ; &&
     x \mapsto (0,x)
\end{align*}
The maps $r$ and $m$ both compose appropriately with the projection map $p$:
\begin{align*}
&p(a,x) = x = p(-a,x) = p(r(a,x)) \\
&p(a_1,x) = p(a_2,x) = x = p(a_1+a_2,x) = p(m((a_1,x),(a_2,x))
\end{align*}
Furthermore, $ps = \Id_X$.  Also
\begin{align*}
m(m((a_1,x),(a_2,x)),(a_3,x)) &= m((a_1+a_2,x),(a_3,x)) \\
&= (a_1+a_2+a_3,x) \\
&= m((a_1,x),(a_2+a_3,x)) \\
&= m((a_1,x),m((a_2,x),(a_3,x))),\\
m(s(x),(a,x)) &= m((0,x),(a,x)) \\
&= (a,x) \\
&= m((a,x),(0,x)) \\
&= m((a,x),s(x)), \\
m((a_1,x),(a_2,x)) &= (a_1+a_2,x) \\
&= (a_2+a_1,x) \\
&= m((a_2,x),(a_1,x)),
\end{align*}
\begin{align*}
\mbox{and}\qquad m(r(a,x),(a,x)) &= m((-a,x),(a,x)) \\
&= (0,x) \\
&= m((a,x),(-a,x)) \\
&= m((a,x),r(a,x)),
\end{align*}
so $T\rmod{A}$ is an Abelian group object in $\categ{Rack}/X$.

Now, given an Abelian group object $p\co R \to X$ in $\categ{Rack}/X$,
with multiplication map $\mu$, inverse map $\nu$, and section $\sigma$,
let $R_x$ be the preimage $p^{-1}(x)$ for each $x \in X$.  Each of the
$R_x$ has a canonical Abelian group structure defined in terms of the maps
$\mu, \nu$, and $\sigma$:  $\sigma(x)$ is the identity in $R_x$, and for
any $u,v \in R_x$ let $u+v := \mu(u,v)$ and $-u := \nu(u)$.  That the preimage
$R_x$ is closed under addition and inversion follows immediately from the
fact that $\mu$ and $\nu$ are rack homomorphisms over $X$.

Next, we define maps
$$\rho_{x,y}\co R_x \to R_{x^y},\ \mbox{given by}\  u \mapsto u^{\sigma(y)},$$
for all $x,y \in X$ and $u \in R_x$.  These are Abelian group
homomorphisms, since $\rho_{x,y}\sigma(x) = \sigma(x)^{\sigma(y)} =
\sigma(x^y)$ (which is the identity in $R_{x^y}$) and, for any $u_1,u_2 \in
R_x$,
\begin{align*}
\rho_{x,y}(u_1+u_2) &= \mu(u_1,u_2)^{\sigma(y)} \\
&= \mu(u_1,u_2)^{\mu(\sigma(y),\sigma(y))} \\
&= \mu(u_1^{\sigma(y)},u_2^{\sigma(y)}) \\
&= \rho_{x,y}(u_1) + \rho_{x,y}(u_2).
\end{align*}
It is also an isomorphism, since exponentiation by a fixed element of a
rack is a bijection.  Furthermore, for any $x,y,z \in X$ and any $u \in
R_x$
\begin{align*}
\rho_{x^y,z}\rho_{x,y}(u) &= u^{\sigma(y)\sigma(z)} \\
&= u^{\sigma(z)\sigma(y)^{\sigma(z)}} \\
&= u^{\sigma(z)\sigma(y^z)} \\
&= \rho_{x^z,y^z}\rho_{x,z}(u).
\end{align*}
Now we define maps
$$\lambda_{y,x}\co R_y \to R_{x^y},\qquad \mbox{given by}\ v \mapsto \sigma(x)^v,$$
for all $x,y \in X$ and $v \in R_y$.  These are also Abelian group
homomorphisms since
$$\lambda_{y,x}\sigma(y) = \sigma(x)^{\sigma(y)} =
\sigma(x^y)$$
(which is the identity in $R_{x^y}$) and, for any $v_1,v_2
\in R_y$,
\begin{align*}
\lambda_{y,x}(v_1+v_2) &= \sigma(x)^{\mu(v_1,v_2)} \\
&= \mu(\sigma(x),\sigma(x))^{\mu(v_1,v_2)} \\
&= \mu(\sigma(x)^{v_1},\sigma(x)^{v_2}) \\
&= \lambda_{y,x}(v_1) + \lambda_{y,x}(v_2).
\end{align*}
Also, for any $x,y,z \in X$, $v \in R_y$ and $w \in R_z$
\begin{align*}
\rho_{x^y,z}\lambda_{y,x}(v) &= \sigma(x)^{v\sigma(z)} \\
&= \sigma(x)^{\sigma(z)v^{\sigma(z)}} \\
&= \sigma(x^z)^{v^{\sigma(z)}} \\
&= \lambda_{y^z,x^z}\rho_{y,z}(v) \\
\mbox{and}\qquad\lambda_{z,x^y}(w) &= \sigma(x^y)^w \\
&= \sigma(x)^{\sigma(y)w} \\
&= \sigma(x)^{w\sigma(y)^w} \\
&= \mu(\sigma(x),\sigma(x))^{\mu(\sigma(z),w)
  \mu(\sigma(y),\sigma(y))^{\mu(w,\sigma(z))}} \\
&= \mu\left(\sigma(x)^{\sigma(z)\sigma(y)^w},
  \sigma(x)^{w\sigma(y)^{\sigma(z)}}\right) \\
&= \sigma(x)^{\sigma(z)\sigma(y)^w} + \sigma(x)^{w\sigma(y)^{\sigma(z)}} \\
&= \sigma(x^z)^{\sigma(y)^w} + \sigma(x)^{w\sigma(y^z)} \\
&= \lambda_{y^z,x^z}\lambda_{z,y}(w) + \rho_{x^z,y^z}\lambda_{z,x}(w).
\end{align*}
Thus an Abelian group object $R \to X$ in $\categ{Rack}/X$ determines a
unique rack module $\rmod{R} = (R,\rho,\lambda)$ over $X$.

For any two such Abelian group objects $p_1\co R_1 \to X$ and $p_2\co R_2
\to X$, together with a rack homomorphism $f_1\co R_1 \to R_2$ over $X$, we
may construct two $X$--modules $\rmod{R}_1$ and $\rmod{R}_2$ as described
above, and an $X$--map $g_1\co\rmod{R}_1\to\rmod{R}_2$ by setting $(g_1)_x(u) =
f_1(u)$ for all $u \in (R_1)_x$ and $x \in X$.  It may be seen that $(g_1)_x:(R_1)_x
\to (R_2)_x$ since $f_1$ is a rack homomorphism over $X$.  It may also be
seen that $g_1$ is a natural transformation of trunk maps $\trunk{T}(X) \to
\trunk{Ab}$ since
\begin{align*}
(g_1)_{x^y}((\rho_1)_{x,y}(u)) &= f_1((\rho_1)_{x,y}(u)) \\
&= f_1(u^{\sigma_1(y)}) \\
&= f_1(u)^{f_1\sigma_1(y)} \\
&= f_1(u)^{\sigma_2(y)} \\
&= (\rho_2)_{x,y}(g_1)_x(u) \\
\mbox{and}\qquad(g_1)_{x^y}(\lambda_1)_{y,x}(v) &= f_1((\lambda_1)_{y,x}(v)) \\
&= f_1(\sigma_1(x)^v) \\
&= f_1\sigma_1(x)^{f_1(v)} \\
&= \sigma_2(x)^{f_1(v)} \\
&= (\lambda_2)_{y,x}(g_1)_x(v)
\end{align*}
for all $u \in R_x, v \in R_y$, and $x,y \in X$.

Given a third Abelian group object $p_3\co R_3 \to X$ together with
another slice morphism $f_2\co R_2 \to R_3$, we may construct another
$X$--module $\rmod{R}_3$ and $X$--map $g_2\co \rmod{R}_2 \to \rmod{R}_3$.
From the composition $f_2f_1$ we may similarly construct a unique $X$--map
$g\co\rmod{R}_1\to\rmod{R}_3$.  Then
$$g_x(u) = (f_2f_1)(u) = (g_2)_x(f_1(u)) = (g_2)_x(g_1)_x(u).$$
Hence this construction determines a functor
$S\co\categ{Ab}(\categ{Rack}/X)\to\categ{RMod}_X$, which is the inverse
of the functor $T\co\categ{RMod}_X\to\categ{Ab}(\categ{Rack}/X)$
described earlier.
\end{proof}

\begin{theorem}
\label{cor:rmod-abelian}
\label{thm:rmod-abelian}
The category $\categ{RMod}_X$ is Abelian.
\end{theorem}
\begin{proof}
The category $\categ{RMod}_X$ is additive, as for any $X$--modules
$\rmod{A}$ and $\rmod{B}$, the set
$\Hom_{\categ{RMod}_X}(\rmod{A},\rmod{B})$ has an Abelian group structure
given by $(f+g)_x(a) = f_x(a) + g_x(a)$ for all
$f,g\co\rmod{A}\to\rmod{B}$, all $x \in X$ and all $a \in A_x$.
Furthermore, composition of $X$--maps distributes over this addition
operation.  The $X$--module with trivial orbit groups and structure
homomorphisms is the zero object in $\categ{RMod}_X$, and for any two
$X$--modules $\rmod{A} = (A,\alpha,\epsilon)$ and $\rmod{B} =
(B,\beta,\zeta)$, the Cartesian product $\rmod{A}\times\rmod{B}
= (A \times B, \alpha\times\beta, \epsilon\times\zeta)$ is also an $X$--module.

Given an $X$--map $f\co\rmod{B} = (B,\beta,\zeta) \to \rmod{C} =
(C,\gamma,\eta)$ let $\rmod{A} = (A,\alpha,\epsilon)$ such that
$A_x = \{a \in B_x : f_x(a) = 0\}$, with $\alpha_{x,y} =
\beta_{x,y}|_{A_x}$ and $\epsilon_{y,x} = \zeta_{y,x}|_{A_y}$.  Then
$\rmod{A}$ is a submodule of $\rmod{B}$ and the inclusion
$\iota\co\rmod{A}\hookrightarrow\rmod{B}$ is the (categorical) kernel of
$f$.

Now define $\rmod{D} = (D,\delta,\xi)$ where $D_x = C_x/\im f_x$, and
$\delta_{x,y} =\gamma_{x,y}+\im f_x$ and $\xi_{y,x} = \eta_{y,x} + \im
f_y$.  Then $\rmod{D}$ is a quotient of $\rmod{C}$ and the canonical
projection map $\pi\co\rmod{C}\to\rmod{D}$ is the (categorical) cokernel of
$f$.

Let $\mu\co\rmod{H}\to\rmod{K}$ be an $X$--monomorphism.  Then the inclusion
$\iota\co\im\mu\to\rmod{K}$ is a kernel of the quotient map
$\pi\co\rmod{K}\to\rmod{K}/\im\mu$.
Since $\mu$ is injective,
$\mu'\co\rmod{H}\iso\im\mu$ where $\mu'_x(a) = \mu_x(a)$ for all $x \in X$
and $a \in H_x$.  But since kernels are unique up to composition with an
isomorphism, and since $\mu = \iota\mu'$, it follows that $\mu$ is the
kernel of its cokernel, the canonical quotient map $\pi$.

Let $\nu\co\rmod{H}\to\rmod{K}$ be an $X$--epimorphism.  Then the inclusion
map $\iota\co\ker\nu\hookrightarrow\rmod{H}$ is a kernel of $\nu$.
Given another $X$--map $\kappa\co\rmod{H}\to\rmod{L}$ such that
$\kappa\iota = 0$, then $\ker\nu\subseteq\ker\kappa$ so that $\nu(a) =
\nu(b)$ implies that $\kappa(a)=\kappa(b)$.  But since $\nu$ is surjective
we can define an $X$--map $\theta\co\rmod{K}\to\rmod{L}$
by $\theta_x\nu_x(a) = \kappa_x(a)$ for all $a \in H_x$ and $x \in X$.
Then $\theta\nu = \kappa$ and so $\nu$ is a cokernel of $\iota$.

So, every $X$--map has a kernel and a cokernel, every monic $X$--map is
the kernel of its cokernel, and every epic $X$--map is the cokernel of its
kernel, and hence $\categ{RMod}_X$ is an Abelian category.
\end{proof}

These results justify the use of the term `rack module' to describe
the objects under consideration, and show that $\categ{RMod}_X$ is an
appropriate category in which to develop homology theories for racks.
Papers currently in preparation will investigate new homology theories
for racks, based on the derived functor approach of Cartan and Eilenberg
\cite{cartan/eilenberg:homalg} and the cotriple construction of Barr
and Beck \cite{barr/beck:standard-constructions}.

We now introduce a notational convenience which may serve to simplify
matters in future.  Let $X$ be a rack, $\rmod{A} = (A,\phi,\psi)$ an
$X$--module, and $w = y_1 y_2 \ldots y_n$ a word in $\As X$.  Then we may
denote the composition
$$\phi_{x^{y_1 \ldots y_{n-1}},y_n} \phi_{x^{y_1 \ldots y_{n-2}},y_{n-1}}
  \ldots \phi_{x,y_1}$$
by $\phi_{x,w} = \phi_{x,y_1 \ldots y_n}$.  This shorthand is well-defined
as the following lemma shows:

\begin{lemma}
\label{lem:asword}
\label{thm:asword}
If $y_1 \ldots y_n$ and $z_1 \ldots z_m$ are two different representative
words for the same element $w \in \As X$, then the compositions
\begin{align*}
&\phi_{x^{y_1 \ldots y_{n-1}},y_n} \phi_{x^{y_1 \ldots y_{n-2}},y_{n-1}}
  \ldots \phi_{x,y_1} \\
\mbox{and}\qquad &\phi_{x^{z_1 \ldots z_{m-1}},z_m} \phi_{x^{z_1 \ldots
  z_{m-2}},z_{m-1}} \ldots \phi_{x,z_1}
\end{align*}
are equal, for all $x \in X$.  Furthermore, $\phi_{x,1} = \Id_{A_x}$,
where $1$ denotes the identity in $\As X$.
\end{lemma}

\begin{proof}
Let $T\co\categ{RMod}_X\to\categ{Ab}(\categ{Rack}/X)$ be the functor
constructed in the proof of Theorem~\ref{thm:rmod-beck}, and recall that
$R_x = T(\rmod{A})_x$ has an Abelian group structure.  For any $x,y \in
X$, the homomorphism $T(\phi_{x,y})\co R_x \to R_{x^y}$ maps $u \mapsto
u^{\sigma(y)}$, where $\sigma$ is the section of $T\rmod{A}$.
Then for any $u \in R_x$
\begin{align*}
T(\phi_{x^{y_1 \ldots y_{n-1}},y_n} \phi_{x^{y_1 \ldots y_{n-2}},y_{n-1}}
  \ldots &\phi_{x,y_1})(u) \\
&= u^{\sigma(y_1)\ldots\sigma(y_n)} \\
&= u^{\sigma(y_1 \ldots y_n)} \\
&= u^{\sigma(z_1 \ldots z_m)} \\
&= u^{\sigma(z_1)\ldots\sigma(z_m)} \\
&= T(\phi_{x^{z_1 \ldots z_{m-1}},z_m} \phi_{x^{z_1 \ldots z_{m-2}},z_{m-1}}
  \ldots \phi_{x,z_1})(u)
\end{align*}
where the equality in the second and third lines follows from the
functoriality of the associated group.

The final statement follows from the observation
$$T(\phi_{x,1})(u) = u^1 = u = T(\Id_{A_x})(u).$$
Hence this notation is well-defined.
\end{proof}

\subsection{Quandle modules}
We now study the specialisation of rack modules to the subcategory
$\categ{Quandle}$.  A \defn{quandle module} is a rack module $\rmod{A} =
(A,\phi,\psi)$ which satisfies the additional criterion
\begin{equation}
\label{eqn:qmod}
\psi_{x,x}(a) + \phi_{x,x}(a) = a
\end{equation}
for all $a \in A_x$ and $x \in X$.  Where the context is clear, we may
refer to such objects as \defn{$X$--modules}.  There is an obvious notion
of a \defn{homomorphism} (or, in the absence of ambiguity, an
\defn{$X$--map}) of quandle modules, and thus we may form the category
$\categ{QMod}_X$ of quandle modules over $X$.

Similarly to example \ref{exm:andr/grana}, Andruskiewitsch and Gra\~na's
definition of quandle modules coincides with the definition of a
homogeneous quandle module in the sense of the current discussion.

Examples \ref{exm:abgroup}, \ref{exm:alexander}, and \ref{exm:dihedral} of
the previous subsection, are also quandle modules.  Example \ref{exm:asx}
is not, but the variant obtained by setting $\psi_{y,x} = \Id_{A} -
\phi_{x,y}$, for all $x,y \in X$, is.

\begin{example}
\label{exm:andr/grana-qmod}
{\rm For an arbitrary quandle $X$, Andruskiewitsch and
Gra\~na \cite{andr/grana:pointed-hopf} further define a \defn{quandle
$X$--module} to be a rack module (as in example \ref{exm:andr/grana})
which satisfies the additional condition
$$\eta_{x,x} + \tau_{x,x} = \Id_A$$
for all $x \in X$.  This may be seen to be a homogeneous quandle
$X$--module in the context of the current discussion.}
\end{example}

Given a quandle $X$ and a quandle $X$--module $\rmod{A}$, the semidirect
product $\rmod{A} \rtimes X$ has the same definition as before.
\begin{proposition}
\label{thm:semidirect-quandle}
\label{prp:semidirect-quandle}
If $X$ is a quandle and $\rmod{A} = (A,\phi,\psi)$ a quandle module over
$X$, the semidirect product $\rmod{A} \rtimes X$ is a quandle.
\end{proposition}
\begin{proof}
By proposition \ref{prp:semidirect-rack}, $\rmod{A} \rtimes X$ is a rack,
so we need only verify the quandle axiom.  For any element $(a,x) \in
\rmod{A} \rtimes X$,
$$(a,x)^{(a,x)} = (\phi_{x,x}(a) + \psi_{x,x}(a),x^x) = (a,x)$$
and so $\rmod{A} \rtimes X$ is a quandle.
\end{proof}

These objects coincide with the Beck modules in the category
$\categ{Quandle}$.

\begin{theorem}
\label{thm:qmod-beck}
For any quandle $X$, there is an equivalence of categories
$$\categ{QMod}_X \iso \categ{Ab}(\categ{Quandle}/X)$$
\end{theorem}
\begin{proof}
As in the proof of Theorem~\ref{thm:rmod-beck}, we identify the quandle
module $\rmod{A} = (A,\phi,\psi)$ with $\rmod{A} \rtimes X \to X$ in the
slice category $\categ{Quandle}/X$.  Proposition
\ref{prp:semidirect-quandle} ensures that this object
is indeed a quandle over $X$, and hence we obtain a well-defined functor
$T\co\categ{QMod}_X\to\categ{Ab}(\categ{Quandle}/X)$.

Conversely, suppose that $R\to X$ is an Abelian group object in
$\categ{Quandle}/X$, with multiplication map $\mu$, inverse map $\nu$, and
section $\sigma$.  As before, we may construct a rack module $\rmod{R} =
(R,\rho,\lambda)$ over $X$.  It remains only to show that this module
satisfies the additional criterion (\ref{eqn:qmod}) for it to be a
quandle module over $X$.  But
\begin{align*}
\lambda_{x,x}(a) + \rho_{x,x}(a) &= \mu(\sigma(x)^a,a^{\sigma(x)}) \\
&= \mu(\sigma(x),a)^{\mu(a,\sigma(x))} \\
&= \mu(\sigma(x),a)^{\mu(\sigma(x),a)} \\
&= \mu(\sigma(x)^{\sigma(x)},a^a) = a
\end{align*}
and so $\rmod{R}$ is indeed a quandle $X$--module.
\end{proof}
\begin{theorem}
The category $\categ{QMod}_X$ is Abelian.
\end{theorem}
\begin{proof}
This proof is exactly the same as the proof of
Theorem~\ref{cor:rmod-abelian}.
\end{proof}
Analogously to the previous subsection, we may conclude that our use of
the term `quandle module' is justified, and that the category
$\categ{QMod}_X$ is a suitable environment in which to study the homology
and cohomology of quandles.

\section{Abelian extensions}
Having characterised suitable module categories, we may now study
extensions of racks and quandles by these objects.  Rack extensions have
been studied before, in particular by Ryder \cite{ryder:thesis} under the
name `expansions'; the constructs which she dubs `extensions' are in some
sense racks formed by disjoint unions, whereby the original rack becomes a
subrack of the `extended' rack.  Ryder's notion of rack expansions is
somewhat more general than the extensions studied here, as she
investigates arbitrary congruences (equivalently, rack epimorphisms onto a
quotient rack) whereas we will only examine certain classes of such
objects.

\subsection{Abelian extensions of racks}

An \defn{extension} of a rack $X$ by an $X$--module $\rmod{A}
= (A,\phi,\psi)$ consists of a rack $E$ together with an
epimorphism $f\co E \epic X$ inducing a partition $E = \bigcup_{x \in X}
E_x$ (where $E_x$ is the preimage $f^{-1}(x)$), and for each $x \in X$ a
left $A_x$--action on $E_x$ satisfying the following three conditions:
\begin{enumerate}
\item[(X1)] The $A_x$--action on $E_x$ is simply transitive, which is to say
that for any $u,v \in E_x$ there is a unique $a \in A_x$ such that $a
\cdot u = v$.
\item[(X2)] For any $u \in E_x$, $a \in A_x$, and $v \in E_y$, $(a \cdot
u)^v = \phi_{x,y}(a) \cdot (u^v)$.
\item[(X3)] For any $u \in E_y$, $b \in A_y$, and $v \in E_y$, $u^{(b \cdot
v)} = \psi_{y,x}(b) \cdot (u^v)$.
\end{enumerate}

Two extensions $f_1\co E_1 \epic X$ and $f_2\co E_2 \epic X$ by the same
$X$--module $\rmod{A}$ are \defn{equivalent} if there exists a rack
isomorphism (an \defn{equivalence}) $\theta\co E_1 \to E_2$ which respects
the projection maps and the group actions:
\begin{enumerate}
\item[(E1)] $f_2\theta(u) = f_1(u)$ for all $u \in E_1$
\item[(E2)] $\theta(a \cdot u) = \theta(a) \cdot u$ for all $u \in E_x$,
$a \in A_x$ and $x \in X$.
\end{enumerate}

Let $f\co E \epic X$ be an extension of $X$ by $\rmod{A}$.  Then a
\defn{section} of $E$ is a function (not necessarily a rack homomorphism)
$s\co X \to E$ such that $fs = \Id_X$.  Since the $A_x$ act simply
transitively on the $E_x$, there is a unique $x \in X$ and a unique $a \in
A_x$ such that a given element $u \in E_x$ can be written as $u = a \cdot
s(x)$.  Since $f$ is a homomorphism, it follows that $s(x)^{s(y)} \in
E_{x^y}$ and so there is a unique $\sigma_{x,y} \in A_{x^y}$ such that
$s(x)^{s(y)} = \sigma_{x,y} \cdot s(x^y) $.  The set $\sigma = \{
\sigma_{x,y} : x,y \in X \}$ is the \defn{factor set} of the extension $E$
\defn{relative to} the section $s$, and may be regarded as an obstruction
to $s$ being a rack homomorphism.

It follows that, for all $x,y \in X$, $a \in A_x$, and $b \in A_y$
\begin{align*}
(a \cdot s(x))^{(b \cdot s(y))} &= \phi_{x,y}(a) \cdot s(x)^{(b \cdot s(y))} \\
&= (\psi_{y,x}(b) + \phi_{x,y}(a)) \cdot s(x)^{s(y)} \\
&= (\psi_{y,x}(b) + \phi_{x,y}(a) + \sigma_{x,y}) \cdot s(x^y) 
\end{align*}
Thus the rack structure on $E$ is determined completely by the factor set
$\sigma$.  The next result gives necessary and sufficient conditions
on factor sets of arbitrary rack extensions.
\begin{proposition}
\label{thm:rack-ext1}
\label{prp:rack-ext1}
Let $X$ be a rack, and $\rmod{A} = (A,\phi,\psi)$ be an $X$--module.  Let
$\sigma = \{ \sigma_{x,y} \in A_{x^y} : x,y \in X \}$ be a collection of
group elements.  Let $E[\rmod{A},\sigma]$ be the set $\{ (a,x) : a \in
A_x, x \in X \}$ with rack operation
$$(a,x)^{(b,y)} = (\phi_{x,y}(a) + \sigma_{x,y} + \psi_{y,x}(b), x^y)$$
for all $a \in A_x$, $b \in A_y$, and $x,y \in X$.

Then $E[\rmod{A},\sigma]$ is an extension of $X$ by $\rmod{A}$ with factor
set $\sigma$ if
\begin{equation}
\label{eqn:rack-ext1}
\sigma_{x^y,z} + \phi_{x^y,z}(\sigma_{x,y}) = \phi_{x^z,y^z}(\sigma_{x,z})
  + \sigma_{x^z,y^z} + \psi_{y^z,x^z}(\sigma_{y,z})
\end{equation}
for all $x,y,z \in X$.
Conversely, if $E$ is an extension of $X$ by $\rmod{A}$ with factor set
$\sigma$ then (\ref{eqn:rack-ext1}) holds, and $E$ is equivalent to
$E[\rmod{A},\sigma]$.
\end{proposition}
\begin{proof}
To prove the first part, we require that $E[\rmod{A},\sigma]$ satisfy the
rack axioms.  Given $(a,x), (b,y) \in E[\rmod{A},\sigma]$, there is a
unique $(c,z) \in E[\rmod{A},\sigma]$ such that $(c,z)^{(b,y)} = (a,x)$,
given by
$$(c,z) = (\phi_{x,y}^{-1}(a - \sigma_{x,y} - \psi_{y,x}(b)),x^{\overline{y}})$$
Also, for any $(a,x), (b,y), (c,z) \in  E[\rmod{A},\sigma]$,
\begin{align*}
(a,x)^{(b,y)(c,z)}
  &= (\phi_{x,y}(a) + \sigma_{x,y} + \psi_{y,x}(b),x^y)^{(c,z)} \\
&= (\phi_{x^y,z}\phi_{x,y}(a) + \phi_{x^y,z}(\sigma_{x,y}) +
  \phi_{x^y,z}\psi_{y,x}(b) + \sigma_{x^y,z} + \psi_{z,x^y}(c), x^{yz})
\end{align*}
and
\begin{align*}
(a,x)^{(c,z)(b,y)^{(c,z)}}
  &= (\phi_{x,z}(a) + \sigma_{x,z} + \psi_{z,x}(c),x^z)^{(\phi_{y,z}(b) +
  \sigma_{y,z} + \psi_{z,y}(c),y^z)} \\
&= (\phi_{x^z,y^z}\phi_{x,z}(a) + \phi_{x^z,y^z}(\sigma_{x,z}) +
  \phi_{x^z,y^z}\psi_{z,x}(c) + \sigma_{x^z,y^z} \\
&\qquad+ \psi_{y^z,x^z}\phi_{y,z}(b) + \psi_{y^z,x^z}(\sigma_{y,z}) +
  \psi_{y^z,x^z}\psi_{z,y}(c), x^{zy^z})
\end{align*}
are equal if (\ref{eqn:rack-ext1}) holds, and so $E[\rmod{A},\sigma]$ is a
rack.

Now define $f\co E[\rmod{A},\sigma] \epic X$ to be projection onto the
second coordinate, and let $A_x$ act on $E[\rmod{A},\sigma]_x = f^{-1}(x)$
by $a_1 \cdot (a_2,x) := (a_1 + a_2,x)$ for each $a_1,a_2 \in A_x$ and all $x
\in X$.  These actions are simply transitive and satisfy the requirements
\begin{align*}
(a_1 \cdot (a_2,x))^{(b,y)} = (a_1+a_2, x)^{(b,y)}
&= (\phi_{x,y}(a_1+a_2) + \sigma_{x,y} + \psi_{y,x}(b),x^y) \\
&= (\phi_{x,y}(a_1) + \phi_{x,y}(a_2) + \sigma_{x,y} + \psi_{y,x}(b),x^y) \\
&= \phi_{x,y}(a_1) \cdot (\phi_{x,y}(a_2) + \sigma_{x,y} +
  \psi_{y,x}(b),x^y) \\
&= \phi_{x,y}(a_1) \cdot (a_2,x)^{(b,y)}\\
\mbox{and}\qquad(a,x)^{b_1 \cdot (b_2,y)} &= (a,x)^{(b_1+b_2,y)} \\
&= (\phi_{x,y}(a) + \sigma_{x,y} + \psi_{y,x}(b_1+b_2),x^y) \\
&= (\phi_{x,y}(a) + \sigma_{x,y} + \psi_{y,x}(b_1) + \psi_{y,x}(b_2),x^y) \\
&= \psi_{y,x}(b_1) \cdot (\phi_{x,y}(a) + \sigma_{x,y} + \psi_{y,x}(b_2),x^y) \\
&= \psi_{y,x}(b_1) \cdot (a,x)^{(b_2,y)}
\end{align*}
so $E[\rmod{A},\sigma]$ is an extension of $X$ by $\rmod{A}$.  Now
define $s\co X \epic E[\rmod{A},\sigma]$ by $s(x) = (0,x)$ for all $x \in
X$.  This is clearly a section of this extension.  Also,
$$s(x)^{s(y)} = (0,x)^{(0,y)} = (\sigma_{x,y},x^y) =
  \sigma_{x,y} \cdot s(x^y)$$
so $\sigma$ is the factor set of this extension relative to the
section $s$.

Conversely, let $f\co E \epic X$ be an extension of $X$ by a given
$X$--module $\rmod{A}$, with factor set $\sigma$ relative to some
extension $s\co X \to E$.  By the simple transitivity of the $A_x$--action
on the $E_x = f^{-1}(x)$, the map $\theta\co (a,x) \mapsto a \cdot s(x)$
is an isomorphism $E[\rmod{A},\sigma] \iso E$.  Since $E$ is a rack, the
earlier part of the proof shows that (\ref{eqn:rack-ext1}) holds, and so
$E[\rmod{A},\sigma]$ is another extension of $X$ by $\rmod{A}$.
Furthermore, $\theta$ respects the projection maps onto $X$, and
$$\theta(a_1 \cdot (a_2,x)) = \theta(a_1+a_2,x) = (a_1+a_2) \cdot s(x)
  = a_1 \cdot (a_2 \cdot s(x)) = a_1 \cdot \theta(a_2,x)$$
so $\theta$ is an equivalence of extensions.
\end{proof}

Andruskiewitsch and Gra\~na \cite{andr/grana:pointed-hopf} introduce the
notion of an extension by a \defn{dynamical cocycle}.  Given an arbitrary
rack $X$ and a non-empty set $S$, we select a function $\alpha\co X \times
X \to \Hom_\categ{Set}(S \times S,S)$ (which determines, for each
ordered pair $x,y \in X$, a function $\alpha_{x,y}\co S \times S \to
S$) satisfying the criteria
\begin{enumerate}
\item $\alpha_{x,y}(s,-)$ is a bijection on $S$
\item $\alpha_{x^y,z}(s,\alpha_{x,y}(t,u)) =
  \alpha_{x^z,y^z}(\alpha_{x,z}(s,t),\alpha_{x,y}(s,u))$
\end{enumerate}
for all $x,y,z \in X$ and $s,t,u \in S$.  Then we may define a rack
structure on the set $X \times S$ by defining $(x,s)^{(y,t)} =
(x^y,\alpha_{x,y}(s,t))$.  This rack, denoted $X \times_\alpha S$, is the
extension of $X$ by $\alpha$.  In the case where $S$ is an Abelian group, and
$\alpha_{x,y}(s,t) = \phi_{x,y}(s) + \sigma_{x,y} + \psi_{y,s}(t)$ for
some suitably-chosen Abelian group homomorphisms $\phi_{x,y},\psi_{y,x}\co S \to S$, and
family $\sigma = \{ \sigma_{x,y} \in S : x,y \in X \}$ of elements of $S$,
then this is equivalent to the construction $E[\rmod{A},\sigma]$ just
discussed, for a homogeneous $X$--module $\rmod{A} = (A,\phi,\psi)$.

\begin{proposition}
\label{thm:rack-ext2}
\label{prp:rack-ext2}
Let $\sigma$ and $\tau$ be factor sets corresponding to extensions of a
rack $X$ by an $X$--module $\rmod{A}$.  Then the following are equivalent:
\begin{enumerate}
\item $E[\rmod{A},\sigma]$ and $E[\rmod{A},\tau]$ are equivalent
  extensions of $X$ by $\rmod{A}$
\item there exists a family $\upsilon = \{ \upsilon_x \in A_x : x \in X \}$
  such that
  \begin{equation}
  \label{eqn:rack-ext2}
  \tau_{x,y} = \sigma_{x,y}
    + \phi_{x,y}(\upsilon_x) + \psi_{y,x}(\upsilon_y) - \upsilon_{x^y}
  \end{equation}
  for $x,y \in X$.
\item $\sigma$ and $\tau$ are factor sets of the same extension of $X$ by
  $\rmod{A}$, relative to different sections.
\end{enumerate}
\end{proposition}
\begin{proof}
Let $\theta\co E[\rmod{A},\tau] \iso E[\rmod{A},\sigma]$ be the
hypothesised equivalence.  Then it follows that $\theta(0,x) =
(\upsilon_x,x)$ for some $\upsilon_x \in A_x$ and, furthermore,
$$\theta(a,x) = \theta(a \cdot (0,x)) = a \cdot \theta(0,x)
  = a \cdot (\upsilon_x,x) = (a + \upsilon_x,x)$$
for all $a \in A_x$, since $\theta$ preserves the $A_x$--actions.
Then
\begin{align*}
&\theta\bigl((a,x)^{(b,y)}\bigr) = (\phi_{x,y}(a) + \psi_{y,x}(b) +
  \tau_{x,y} + \upsilon_{x^y}, x^y) \\
\mbox{and}\qquad&\theta(a,x)^{\theta(b,y)} = (a + \upsilon_x,x)^{(b+\upsilon_y,y)} =
  (\phi_{x,y}(a+\upsilon_x) + \psi_{y,x}(b+\upsilon_y) + \sigma_{x,y},x^y)
\end{align*}
which are equal since $\theta$ is a rack homomorphism, and so
(\ref{eqn:rack-ext2}) holds.  This argument is reversible, showing the
equivalence of the first two statements.

Now, given such an equivalence $\theta$, define a section $s\co X \to
E[\rmod{A},\tau]$ by $x \mapsto (\upsilon_x,x)$.  Then the above
argument also shows that
\begin{align*}
&s(x)^{b \cdot s(y)} = (\upsilon_x,x)^{(\upsilon_y+b,y)} =
  (\sigma_{x,y} + \psi_{y,x}(b)) \cdot s(x^y) \\
\mbox{and}\qquad &(a \cdot s(x))^{s(y)} = (\upsilon_x+a,x)^{(\upsilon_y,y)} =
  (\sigma_{x,y} + \phi_{x,y}(a)) \cdot s(x^y)
\end{align*}
so $\sigma$ is the factor set of $E[\rmod{A},\tau]$ relative to the
section $s$.  This property holds for any extension equivalent to
$E[\rmod{A},\tau]$.  Conversely, if $\sigma$ and $\tau$ are factor sets of
some extension $E$ of $X$ by $\rmod{A}$ relative to different sections
$s,t\co X \to E$ then $s(x) = \upsilon_x \cdot t(x)$ for some $\upsilon_x
\in A_x$, and so the first and third conditions are equivalent.
\end{proof}
The following corollary justifies the earlier assertion that the factor
set is in some sense the obstruction to a section being a rack
homomorphism.
\begin{corollary}
\label{thm:rack-ext3}
\label{cor:rack-ext3}
For an extension $f\co E \epic X$ by an $X$--module $\rmod{A} =
(A,\phi,\psi)$, the following statements are equivalent:
\begin{enumerate}
\item There exists a rack homomorphism $s\co X \to E$ such that $fs = \Id_X$
\item Relative to some section, the factor set of $E \epic X$ is trivial
\item Relative to any section there exists, for the factor set $\sigma$ of
$E \epic X$, a family $\upsilon = \{\upsilon_x \in A_x : x \in X \}$ such
that for all $x,y \in X$
\begin{equation}
\label{eqn:rack-ext3}
\sigma_{x,y} = \phi_{x,y}(\upsilon_x) - \upsilon_{x^y} + \psi_{y,x}(\upsilon_y)
\end{equation}
\end{enumerate}
\end{corollary}
Extensions of this type are said to be \defn{split}.  We are now able to
classify rack extensions:
\begin{theorem}
\label{thm:rack-ext4}
Let $X$ be a rack and $\rmod{A} = (A,\phi,\psi)$ an $X$--module.  Then
there is an Abelian group $\Ext(X,\rmod{A})$ whose elements are in
bijective correspondence with extensions of $X$ by $\rmod{A}$.
\end{theorem}
\begin{proof}
Let the set $Z(X,\rmod{A})$ consist of extensions of $X$ by $\rmod{A}$.
As shown above, these are determined by factor sets $\sigma$ satisfying
(\ref{eqn:rack-ext1}).  Defining an addition operation by
$(\sigma+\tau)_{x,y} := \sigma_{x,y} + \tau_{x,y}$ gives this an Abelian
group structure with the trivial factor set as identity.  A routine
calculation confirms that the set
$B(X,\rmod{A})$ of split
extensions (equivalently, factor sets satisfying (\ref{eqn:rack-ext3}))
forms an Abelian subgroup of $Z(X,\rmod{A})$, and so we may define
$\Ext(X,\rmod{A}) := Z(X,\rmod{A})/B(X,\rmod{A})$.
\end{proof}

In the case where $\rmod{A}$ is a trivial homogeneous $X$--module
(equivalently, an Abelian group $A$) the group $\Ext(X,\rmod{A})$
coincides with $H^2(BX;A)$, the
second cohomology group of the rack space of $X$ as defined by Fenn, Rourke and
Sanderson \cite{fenn/rourke/sanderson:trunks}.

\subsection{Abelian extensions of quandles}
We now turn our attention to the case where $X$ is a quandle.  Extensions
of $X$ by a quandle $X$--module $\rmod{A}$ and their
corresponding factor sets are defined in an analogous manner.
\begin{proposition}
\label{thm:quandle-ext1}
\label{prp:quandle-ext1}
Let $X$ be a quandle and $\rmod{A} = (A,\phi,\psi)$ be a quandle
module over $X$.  Then extensions $f\co E \epic X$ such that $E$ is also a
quandle are in bijective correspondence with factor sets $\sigma$
satisfying hypothesis (\ref{eqn:rack-ext1}) of proposition
\ref{thm:rack-ext1} together with the additional criterion
\begin{equation}
\label{eqn:quandle-ext1}
\sigma_{x,x} = 0
\end{equation}
for all $x \in X$.
\end{proposition}
\begin{proof}
Following the reasoning of proposition \ref{thm:rack-ext1}, for $E$ to be
a quandle is equivalent to the requirement that
$$(a,x)^{(a,x)} = (\phi_{x,x}(a) + \sigma_{x,x} + \psi_{x,x}(a),x^x)
  = (a,x)$$
for all $x \in X$ and $a \in A_x$.  Since $\rmod{A}$ is a quandle module,
this is equivalent to the requirement that (\ref{eqn:quandle-ext1}) holds.
\end{proof}
We may now classify quandle extensions of $X$ by $\rmod{A}$:
\begin{theorem}
\label{thm:quandle-ext2}
For any quandle $X$ and quandle $X$--module $\rmod{A}$, there is an
Abelian group $\Ext_Q(X,\rmod{A})$ whose elements are in bijective
correspondence with quandle extensions of $X$ by $\rmod{A}$.
\end{theorem}
\begin{proof}
We proceed similarly to the proof of Theorem~\ref{thm:rack-ext4}.  Let
$Z_Q(X,\rmod{A})$ be the subgroup of $Z(X,\rmod{A})$ consisting of
factor sets satisfying the criterion (\ref{eqn:quandle-ext1}), and let
$B_Q(X,\rmod{A}) = B(X,\rmod{A})$.  Then we define $\Ext_Q(X,\rmod{A}) =
Z_Q(X,\rmod{A})/B_Q(X,\rmod{A})$.
\end{proof}
In the case where $\rmod{A}$ is trivial homogeneous (and hence equivalent
to an Abelian group $A$), extensions of $X$ by $\rmod{A}$ correspond to
Abelian quandle extensions, in the sense of Carter, Saito and
Kamada \cite{carter/kamada/saito:diag} and so $\Ext_Q(X,\rmod{A}) =
H^2_Q(X;A)$.

If the module $\rmod{A}$ is a homogeneous Alexander module as defined
in example \ref{exm:alexander}, then extensions of $X$ by $\rmod{A}$
are exactly the twisted quandle extensions described by Carter,
Saito and Elhamdadi \cite{carter/elhamdadi/saito:twisted}, and so
$\Ext_Q(X,\rmod{A}) = H^2_{TQ}(X;A)$.


\begin{thebibliography}{99}

\bibitem{andr/grana:pointed-hopf}
{Nicol\'as Andruskiewitsch}, {Mat\'ias Gra\~na},
\emph{From racks to pointed Hopf algebras},
%{\tt arXiv:math.QA/0202084}
Advances in Mathematics \textbf{178} (2003) 177--243
%  \MR{1994219}

\bibitem{barr/beck:standard-constructions}
{Michael Barr}, {Jonathan Beck},
\emph{Homology and standard constructions},
from: ``Seminar on Triples and Categorical Homology Theory'', volume 80 of
Lecture Notes in Mathematics, Springer--Verlag (1969) 245--335
%  \MR{0258917}

\bibitem{beck:thesis}
{Jonathan Beck},
\emph{Triples, algebras and cohomology},
PhD thesis, Columbia University (1967).  Republished as: Reprints in
Theory and Applications of Categories \textbf{2} (2003) 1--59
%  \MR{1987896}

\bibitem{brieskorn:automorphic}
{E Brieskorn},
\emph{Automorphic sets and singularities},
Contemporary Mathematics \textbf{78} (1988) 45--115
%  \MR{0975077}

\bibitem{cartan/eilenberg:homalg}
{Henri Cartan}, {Samuel Eilenberg},
\emph{Homological Algebra},
Princeton University Press (1999)
%  \MR{1731415}

\bibitem{carter/elhamdadi/saito:twisted}
{J Scott Carter}, {Mohamed Elhamdadi}, {Masahico Saito},
\emph{Twisted quandle homology theory and cocycle knot invariants},
Algebraic and Geometric Topology \textbf{2} (2002) 95--135
%  \MR{1885217}

\bibitem{carter/kamada/saito:diag}
{J Scott Carter}, {Seiichi Kamada}, {Masahico Saito},
\emph{Diagrammatic computations for quandles and cocycle knot invariants},
%{\tt arXiv:math.GT/0102092}
from: ``Diagrammatic morphisms and applications (San Francisco, CA,
2000)'', Contemporary Mathematics \textbf{318} (2003) 51--74
%  \MR{1973510}

\bibitem{conway/wraith:wracks}
{John Conway}, {Gavin Wraith},
unpublished correspondence (1959)

\bibitem{etingof/grana:orc}
{Pavel Etingof}, {Mat\'ias Gra\~na},
\emph{On rack cohomology},
Journal of Pure and Applied Algebra \textbf{177} (2003) 49--59
%  \MR{1948837}

\bibitem{fenn/rourke:racks-links}
{Roger Fenn}, {Colin Rourke},
\emph{Racks and links in codimension 2},
Journal of Knot Theory and its Ramifications \textbf{1} (1992) 343--406
%  \MR{1194995}

\bibitem{fenn/rourke/sanderson:trunks}
{Roger Fenn}, {Colin Rourke}, {Brian Sanderson},
\emph{Trunks and classifying spaces},
Applied Categorical Structures \textbf{3} (1995) 321--356
%  \MR{1364012}

\bibitem{jackson:thesis}
{Nicholas Jackson},
\emph{Homological algebra of racks and quandles},
PhD thesis, Mathematics Institute, University of Warwick (2004)

\bibitem{joyce:knot-quandle}
{David Joyce},
\emph{A classifying invariant of knots: the knot quandle},
Journal of Pure and Applied Algebra \textbf{23} (1982) 37--65
%  \MR{0638121}

\bibitem{ryder:thesis}
{Hayley Ryder},
\emph{The structure of racks},
PhD thesis, Mathematics Institute, University of Warwick (1993)

\end{thebibliography}
\end{document}